\input amstex 
\documentstyle{amsppt} 
\magnification=\magstep1
\vcorrection{-0.8cm} 
\NoBlackBoxes 
\nologo

\def\({\left(} 
\def\){\right)}

\def\g{\frak g} 
\def\h{\frak h}
\def\n{\frak n} 
\def\bb{\frak b} 
\def\D{\Delta} 
\def\l{\lambda}
\def\Dp{\Delta^+} 
\def\Dm{\Delta^-} 
\def\Da{\widehat\Delta}
\def\Dap{\widehat\Delta^+} 
\def\d{\delta}

\def\a{\alpha} 
\def\b{\beta} 
 
\def\I{\Cal I} 
\def\l{\lambda}
\def\d{\delta} 
 
\def\th{\theta} 
 
\def\i{{\frak i}}

\def\p{\Phi} 
\def\nat{\Bbb N} 
\def\ganz{\Bbb Z} 
 
\def\real{\Bbb R}

\TagsOnRight

\def\l{\lambda} 

\def\({\left(} \def\){\right)}

\topmatter 

\author  
Paola Cellini\\ Paolo Papi
\endauthor 
\dedicatory To Claudio Procesi, on the occasion of his 60$^{th}$ birthday\enddedicatory
\address 
Paola Cellini \vskip 0.pt Dipartimento di Scienze
 \vskip 0.pt Universit\`a ``G. d'Annunzio" \vskip 0.pt Viale Pindaro 42
\vskip 0.pt 65127 Pescara --- ITALY \vskip 0.pt e-mail:{\rm \
cellini\@gotham.sci.unich.it } 
\endaddress 

\bigskip 
\bigskip 
\address Paolo Papi
\vskip 0.pt Dipartimento di Matematica Istituto G. Castelnuovo \vskip 0.pt
Universit\`a di Roma "La Sapienza" \vskip 0.pt Piazzale Aldo Moro 5 \vskip 0.pt
00185 Rome --- ITALY \vskip 0.pt e-mail:{\rm \ papi\@mat.uniroma1.it}
\endaddress

\leftheadtext 
{Paola Cellini Paolo Papi} 
\rightheadtext 
{$ad$-nilpotent  ideals II} 

\title $ad$-Nilpotent ideals of a Borel subalgebra II\endtitle

\abstract We provide an explicit bijection between the $ad$-nilpotent
ideals of a Borel subalgebra of a simple Lie algebra $\g$ and the orbits
of $\check Q/(h+1)\check Q$ under the Weyl group ($\check Q$ being the
coroot lattice and $h$ the Coxeter number of $\g$). From this result we
deduce in a uniform way a counting formula for the $ad$-nilpotent ideals.
\endabstract



\endtopmatter 

\document

\heading 
\S1 Introduction
\endheading

Let $\g$ be a complex simple Lie algebra of rank $n$. Let $\bb\subset\g$
be a fixed Borel subalgebra, with Cartan component $\h$, and let $\Dp$ be
the positive system of the root system $\D$ of $\g$ corresponding to the
previous choice. For each $\a\in \Dp$ let ${\frak g}_\a$ denote the root
space of $\g$ relative to $\a$, and $\n = \bigoplus \limits_{\a\in
\Dp}{\frak g}_\a$, so that $\bb=\h\oplus \n$. \par
In this paper we continue the analysis, started in \cite{CP1}, of the set
$\I$ of $ad$-nilpotent ideals of $\bb$, i.e. the
ideals of $\bb$ consisting of $ad$-nilpotent elements.  These ideals are
precisely the ideals of $\bb$ which are contained in $\n$; in particular
the abelian ideals of $\bb$ are $ad$-nilpotent.
\par
The abelian ideals of Borel subalgebras were studied by
Kostant in \cite{Ko1}, \cite{Ko2} in connection with representation theory of
compact semisimple Lie groups, and very recently by Panyushev and R\"ohrle, \cite{PR}, in connection with the theory of spherical orbits. In particular, in \cite{Ko2}  Kostant detailed the following 
unpublished theorem of D. Peterson: the abelian ideals of $\bb$ are $2^n$ in number, 
independently of the type of $\g$. In fact, Peterson gives a bijection between the abelian ideals of $\bb$ and a certain set
 of elements of the affine Weyl group $\widehat W$ of $\g$. This leads to look for 
a similar result for all $ad$-nilpotent ideals. 
In \cite{CP1} we showed how to associate to any ideal $\i$ in $\I$ a uniquely determined 
element $w_\i$ in $\widehat W$, and we characterized the set $\{w_\i\mid \i\in \I\}$
inside $\widehat W$. In this paper we develop our previous result and prove the following  Theorem. \par
\par Let $W$
denote the Weyl group of $\g$, and  $Q,\,\check Q$ be the root and coroot
lattices, respectively. Moreover  let $h$ be the Coxeter number of $W$ and 
$e_1,\ldots, e_n$  be the  exponents of $W$ \cite{Hu, Section 3.19}.
  We
consider the natural action of $W$ on $\check Q/(h+1)\check Q$. 
\proclaim {Theorem 1} There exists an explicit bijection between the set of 
$ad$-nilpotent ideals of $\bb$ and the set of  $W$-orbits in $\check Q/(h+1)\check Q$.
In particular the number of the $ad$-nilpotent ideals of $\bb$ is
$$\frac{1}{|W|}\prod_{i=1}^n(h+e_i+1). \tag 1$$
\endproclaim\noindent
The fact that formula (1) counts  $W$-orbits in $\check Q/(h+1)\check Q$ follows from a Theorem of M. Haiman \cite{Ha, Theorem 7.4.4}.\par
As we shall recall in Section 4, the $ad$-nilpotent ideals of 
$\bb$ are naturally in bijection with the antichains of the root poset $(\Dp,\leq)$, 
hence with the {\sl $\oplus$--sign types} of $\check\Delta$, and
with the regions of  the Catalan hyperplane arrangement which 
are contained in the fundamental chamber of $W$. 
So our result affords a uniform enumeration for all these objects. In particular, it answers the question raised in \cite{S, Remark 3.7} regarding the (uniform) enumeration of sign types. 
\par
Formula (1) for the $ad$-nilpotent ideals already appears at the end of the Introduction of \cite{KOP}, where the authors asked for a proof  of it avoiding  case by case inspection.
Our main Theorem also improves  the known results  on the Catalan arrangement
extending to any root system the enumeration formula (1),  which 
was proved by Athanasiadis for the classical systems \cite{At1}, \cite{At2}. 
\par
We remark that 
Sommers \cite{So} gives a generalization of formula $(1)$, expressing the Euler
characteristic of the space of partial flags containing a certain regular
semisimple nil-elliptic element $n_t$ in an affine Lie algebra.\smallskip
\par 
Our paper is organized as follows. 
In the next section we fix notation and recall some basic facts about affine Weyl groups. 
In Section 3 we prove Theorem 1 after recalling the results of \cite{CP1} 
which are needed for the proof. In Section 4 we briefly 
recall the known bijections between $ad$-nilpotent ideals, antichains of the root poset
$(\Dp,\leq)$, {\sl $\oplus$--sign types} of $\check\Delta$, and
regions of  the Catalan hyperplane arrangement contained in the fundamental chamber.
In Section 5 we illustrate the bijection of Theorem 1 for the root types $A_2$ and $B_2$.

\medskip 
\heading 
\S 2 Notation and preliminaries
\endheading

\medskip 
Let $\Pi=\{\a_1,\dots,\a_n\}$ be the simple roots of
$\Dp$.  We set $V\equiv \h_{\Bbb R}^*=\bigoplus\limits_{i=1}^n \real \a_i$ and
denote by $(\ ,\ )$ the positive symmetric bilinear form induced on $V$ by the
Killing form. We describe the affine root system associated to $\D$ as follows
\cite{Kac, Chapter 6}. We extend $V$ and its inner product setting $\widehat V=
V \oplus \real \d\oplus \real \l$, $(\d,\d)=(\d,V)=(\l,\l)=(\l,V)=0$, and
$(\d,\l)=1$. We still denote by $(\ ,\ )$ the resulting (non-degenerate)
bilinear form. The affine root system associated to $\D$ is $\Da= \D+\ganz
\d=\{\a+k\d\mid \a\in \D, \ k\in \ganz\}$; remark that the affine roots are
non-isotropic vectors. The set of positive affine roots is $\Dap=(\Dp+\nat
\d)\cup (\Dm+\nat^+\d)$, where $\Dm=-\Dp$. We denote by $\theta$ the highest
root of $\D$ and set $\a_0=-\theta+\d$, $\widehat \Pi=\{\a_0,\dots,\a_n\}$. For
each $\a\in\Dap$ we denote by $s_\a$ the corresponding reflection of $\widehat
V$; the affine Weyl group associated to $\D$ is the group $\widehat W$ generated
by $\{s_\a\mid \a\in \Dap\}$. Note that $w(\delta)=\delta$ for each $w\in
\widehat W$. 
\par 
$\widehat W$ is a semidirect product $T\rtimes W$, 
where $T=\{t_\tau \mid \tau\in\check Q\} \cong \check Q$ is the subgroup of 
{\it translations}, and the action of $W$ on $T$ is $vt_\tau v^{-1}=t_{v(\tau)}$. 
The action of $t_\tau$ on $V\oplus \real \delta$, in particular on the roots, is
given by $t_\tau(x)=x-(x, \tau) \delta$, for each $x\in V\oplus \real \delta$.
(See \cite{Kac} for the general definition of $t_\tau$ on $\widehat V$). 
\par 
Consider the $\widehat W$-invariant
affine subspace $E=\{x\in V\mid (x,\d)=1\}=V\oplus\Bbb R\d+\l$. Let $\pi:E \to
V+\lambda $ be the natural projection. For $w\in \widehat W$ we set $\overline
w=\pi\circ  w_{|V+\lambda}$. Then the map $w\mapsto \overline w$ is injective.
We identify $V+\l$ with $V$ through the natural projection. In this way
$w\mapsto \overline w$ induces an isomorphism of $\widehat W$ onto a group 
$W_{af}$ of affine transformations of $V$, which is in fact the usual affine
representation of the affine Weyl group \cite{B, VI, \S 2}. For $v\in W$, $\overline
v$ is simply the restriction of $v$ to $V$, so we omit the bar. For $\tau\in
\check Q$, $\overline t_\tau$ is the true translation by $\tau$. 

\par 
For $k\in \nat^+$, we set 
$$C_k=\{x\in  C_\infty \mid (x,\theta)< k\},$$ 
$$C_\infty=\{x\in V\mid (x,\alpha_i)>0 \text{ for each } i\in \{1,\dots,n\}\}.$$ 
$C_\infty$ is the fundamental chamber of $W$ in $V$, and
$C_1$ is the fundamental alcove of $W_{af}$ in $V$. The closures 
$\overline C_\infty,\,\overline C_1$ are fundamental domains for the actions
on $V$ of $W,\,W_{af}$ respectively.
\par
As usual we denote by $\{\check\omega_1, \dots, \check\omega_n\}$ 
the dual basis of $\{\alpha_1, \dots, \alpha_n\}$ and we set 
$o_i=\check\omega_i/m_i$, where $\theta=m_1\alpha_1+\cdots+ m_n\alpha_n$. 
For each $k\in \nat^+$,
$\overline C_k$($=k\overline C_1$) is the simplex whose vertices are  
$0, ko_1, \dots, ko_n$. 
\par
Let $\overline T^k=\{\overline t_\tau\mid \tau\in k\check Q\}$ and set $W_k=\overline
T^k\rtimes W$. Note that $W_k$ is the affine Weyl group of ${1\over k}\Delta$ and
$C_k$ is its fundamental chamber with respect to ${1\over k}\Pi$. In particular, 
$\check Q\cap\overline C_k$  is a set of representatives of the orbits of 
$\check Q$ under the natural action of $W_k$. 
For $\tau \in \check Q$, consider its orbit $W_k(\tau)$. 
We have $W_k(\tau)=W(\tau)+k\check Q=W(\tau+k\check Q)$,
hence the orbits of  of $\check Q$ under $W_k$ naturally correspond to the
orbits of $\check Q/k\check Q$ under the action of $W$ [Ha]. In fact, 
in order to prove Theorem 1, we shall prove
that $\I$ is in bijection with $\check Q\cap \overline C_{h+1}$.

\heading
\S3 Proof of Theorem 1
\endheading
\medskip
 
In \cite{CP1} we found an explicit  encoding of the elements of $\I$ by means of 
certain elements of $\widehat W$. We briefly recall this result.
By definition any $ad$-nilpotent ideal $\i$ of $\bb$ is
included in $\n$. Such an ideal $\i$  is, in particular, $ad(\h)$-stable, so
there exists $\p_\i\subseteq \Dp$ such that $\i=\bigoplus\limits_{\a\in\p_\i}\g_\a$. 
We set $\p_\i^{1}=\p_\i$ and $\p_\i^{k+1}=(\p_\i^k+\p_\i)\cap\Dp$, for each $k\in \nat^+$, 
so that $\bigoplus\limits_{\a\in\p_\i^k}\g_\a$ equals $\i^{(k)}$, the $k$-th element
of the descending central series of $\i$. 
Then we associate to $\i$ the following set of positive affine roots: 
$L_\i=\bigcup\limits_{k\geq 1}(-\p_\i^k+k\d)$. In \cite{CP1, Proposition 2.4} 
we proved that there exists a (unique) $w_\i\in\widehat W$ such that 
$L_\i=N(w_\i)=\{\a\in\widehat\Delta^+\mid w_\i^{-1}(\a)<0\}$. Then $\i\mapsto w_\i$ is the
required encoding. We remark that $w_\i$ is explicitly  determined by $L_\i$.  
We also gave the following
characterization,   which will play a crucial role in the sequel.

\proclaim
{Proposition 1}\cite{CP1, Proposition 2.12} Let $w\in \widehat  W$. Then
$w=w_{\frak i} $ for some $\frak i\in \Cal I$ if and only if the following
conditions hold:
\roster 
\item"(a)" $w^{-1}(\alpha)>0$ for each $\alpha\in \Pi$;
\item"(b)"  if $w(\alpha)<0$ for some $\alpha\in \widehat  \Pi$, then $w(\alpha)
=\beta -\delta$ for some $\beta\in \Delta^+$.
\endroster 
\endproclaim

\smallskip

For $\alpha\in \Delta^+$ and $l\in \ganz$ set 
$H_{\alpha,l}=\{x\in V\mid (x,\alpha)=l\}$.
We recall that, for $\alpha\in \Delta^+$, $l\in \nat^+$, $m\in \nat$, 
$w\in \widehat W$,
we have $w^{-1}(-\alpha+l \delta)<0$ if and only if $H_{\alpha,l}$ separates 
$ C_1$ and $\overline w(C_1)$, and 
$w^{-1}(\alpha+m \delta)<0$ if and only if $H_{\alpha,-m}$ separates 
$C_1$ and $\overline w(C_1)$.
From Proposition 1 we obtain the following characterization.

\proclaim
{Proposition 2} Let $w\in \widehat  W$, $w=t_\tau v$, $\tau\in \check
Q$, $v\in W$. Set $\beta_i=v(\alpha_i)$ for $i\in\{1,\dots,n\}$. Then
$w=w_{\frak i} $ for some $\frak i\in \Cal I$ if and only if the following
conditions hold:\roster \item"(i)" $\overline w(C_1)\subset C_{\infty}$; \item"(ii)"
$(\tau,\beta_i) \leq 1$ for each $i\in\{1,\dots,n\}$ and $(\tau, v(\theta)) \geq
- 2$.\endroster 
\endproclaim

\demo 
{Proof} It is clear that condition (a) of Proposition 1 is equivalent to
$\overline w(C_1)\subset C_\infty$. 
\par 
Assume $w=w_{\frak i} $ for some $\frak i\in
\Cal I$. Then (i) holds and this implies, in particular, that $\tau\in \overline
C_\infty$. Since $\beta_i\in \Delta$, if $\beta_i<0$  then $(\tau,\beta_i)\leq
0$. So we may assume $\beta_i>0$. We have
$w(\alpha_i)=\beta_i-(\tau,\beta_i)\delta$, hence if $(\tau,\beta_i) > 0$ we
obtain $w(\alpha_i)<0$. Then by Proposition 1 $w(\alpha_i)=\beta_i-\delta$ and
thus $(\tau,\beta_i)=1$. Then we consider $v(\theta)$. We have $w(\alpha_0)=
-v(\theta)+((\tau,v(\theta))+1)\delta$. If $(\tau,v(\theta))< -1$, then
$w(\alpha_0)<0$, hence by Proposition 1 $(\tau,v(\theta))+1=-1$, hence
$(\tau,v(\theta))=-2$.
\par 
Conversely, assume that (i) and (ii) hold.  Then
condition (a) of Proposition 1 holds. Let $1\leq i\leq n$ and  $w(\alpha_i)<0$.
Then either $(\tau,\beta_i)>0$, or $(\tau,\beta_i)=0$ and $\beta_i<0$. The
latter case cannot occur, otherwise, for $x\in C_1$ we would have $(\overline
w(x),\beta_i)=(v(x),\beta_i)=(x,\a_i)>0$, which is impossible, since $\overline w(x)$
belongs to $\overline C_\infty$  and $\beta_i$ is negative. So we have 
$(\tau,\beta_i)>0$,  hence, by assumption $(\tau,\beta_i)=1$, so that
$w(\alpha_i)=\beta_i-\delta$. Finally assume  $w(\alpha_0)<0$. Then either
$(\tau,v(\theta))+1<0$, or $(\tau,v(\theta))+1=0$ and $v(\theta)>0$. As above we
see that the latter case cannot occur, so, by assumption, $(\tau,v(\theta))=-2$.
This implies $v(\theta)<0$  and $w(\alpha_0)= -v(\theta)-\delta$, hence the
claim. 
\qed
\enddemo

Set $$D=\{\tau\in \check Q \mid (\tau,\alpha_i)\leq 1 \text{ for each } i\in
\{1,\dots ,n\} \text{ and } (\tau,\theta)\geq -2\}.$$ Assume $w_\i=t_{\tau_\i} v_\i$ for some $\frak i\in \Cal I$, $\tau_\i\in
\check Q$,
$v_\i\in W$. Then by Proposition 2 we have  $(\tau_\i,\beta_j) \leq 1$ for each $j\in\{1,\dots,n\}$ and
$(\tau_\i, v_\i(\theta)) \geq - 2$, hence  $(v_\i^{-1}(\tau_\i),\alpha_j) \leq 1$ for each
$j\in\{1,\dots,n\}$ and $(v_\i^{-1}(\tau_\i), \theta) \geq - 2$. It follows that
$t_{\tau_\i} v_\i\mapsto v_\i^{-1}(\tau_\i)$ is a map from  $\{w_\frak i \mid \frak i\in \Cal
I\}$ to $D$.

\proclaim 
{Proposition 3}  The map $F: w_\i=t_{\tau_\i} v_\i\mapsto
v_\i^{-1}(\tau_\i)$ is a bijection between\break   $\{w_\frak i \mid \frak i\in \Cal I\}$ and
$D$. 
\endproclaim

\demo
{Proof} Set, for notational simplicity, $w_\frak i=t_{\tau} v$, $w_\frak j=t_{\sigma} u$ for some $\frak i$
and $\frak j$ in $\Cal I$, $\tau,\sigma \in \check Q$ and $v,u\in W$. Assume
$v^{-1}(\tau)=u^{-1}(\sigma)$. Since $\tau, \sigma\in \overline C_\infty $,
which is a fundamental domain for $W$, we have $\tau=\sigma$ and
$v u^{-1}(\tau)=\tau$. Hence $\overline{t_{\tau}} v (C_1)=\overline{t_{\tau}}v u^{-1}u (C_1) =
v u^{-1} \left(\overline{t_{\tau}} u (C_1) \right)=v u^{-1} \left(\overline{t_{\sigma}} u (C_1) \right)\subset
v u^{-1}(C_\infty)$. But $\overline{t_{\tau}} v (C_1)\subset C_\infty$, hence $v u^{-1}=1$. Thus $F$ is
injective. Next let
$\sigma\in D$. We first see that there exists $v\in W$ such that $t_{v(\sigma)}
v(C_1)\subset C_\infty$: simply take the unique $v\in W$ such that
$v(\sigma+C_1)\subset C_\infty$. Now it is immediate that, since $\sigma\in D$,
$t_{v(\sigma)}v$ also satisfies condition (ii) of Proposition 2, hence
$t_{v(\sigma)}v=w_\frak i$ for some $\frak i$ in $\Cal I$. It is obvious that 
$F$ maps $t_{v(\sigma)}v$ to $\sigma$, thus $F$ is surjective. 
\qed
\enddemo
\remark{\bf Remark} In a forthcoming paper \cite{CP2} we provide characterizations for the elements of $D$ corresponding
through
$F$ to abelian ideals and, among them, for those encoding maximal abelian ideals.\endremark

\smallskip
Let $\check P=\ganz \check\omega_1+ \cdots +\ganz \check\omega_n$
be the coweight lattice of $W$. We denote by  $W_{af}'$ the extended  affine Weyl
group, $W_{af}'=\overline T'\rtimes W$, with $\overline T'=\{\overline t_\tau\mid \tau\in \check
P\}$,  $\overline t_\tau$  the translation by $\tau$. 
As usual, we set $f=[W_{af}':W_{af}]=[\check P:\check Q]$.

\proclaim
{Lemma 1} 
Assume that $k$ and $f$  are relatively prime. Then for each 
$w'\in W_{af}'$ there exists $w\in W_{af}$ such that
$w'(C_k)=w(C_k)$. 
\endproclaim

\demo
{Proof} Let $\theta=\sum_{i=1}^n m_i \alpha_i$ and $J=\{i\mid m_i=1\}$. 
By [IM, Sections 1.7 and 1.8], $\{0\}\cup\{\check\omega_j\mid j\in J\}$ 
is a set of representatives of $\check
P/\check Q$. Moreover, for  each $j\in J$,  $C_k=t_{k\check\omega_j}w_0^jw_0(C_k)$, where
$w_0$ is the longest element of $W$ and $w_0^j$ is the longest element in the
maximal parabolic subgroup of $W$ generated by the reflections with respect to the 
$\a_i$ with $i\not=j$. It
suffices to prove the lemma for $w'\in \overline T'$;  let $w'=\overline t_\sigma$ with 
$\sigma\in \check P$. Then we have $w'(C_k)=\overline t_\sigma(C_k)=
\overline t_{\sigma+k\check\omega_j}w_0^jw_0(C_k)$, for each $j\in J$. If $k$ and
$f$ are relatively prime, then $\{0\}\cup\{ k\check\omega_j\mid j\in J\}$ and hence
$\{\sigma\}\cup\{\sigma+k\check\omega_j\mid j\in J\}$ still are sets of representatives of
$\check P/\check Q$. It follows that exactly one element in
$\{\sigma\}\cup\{\sigma+k\check\omega_j\mid j\in J\}$ belongs to $\check Q$, hence one
among $\overline t_\sigma$, $\overline t_{\sigma+k\check\omega_j}w_0^jw_0$, $j\in J$, 
belongs to $W_{af}$.
\qed
\enddemo

\remark
{\bf Remark} A direct check shows that the prime divisors of $f$ also
divide the Coxeter number of $W$. Hence the assumption of Lemma 1 is satisfied
by any integer $k$ relatively prime to $h$. 
\endremark

\smallskip
\proclaim{Theorem 1} $\Cal I$ is in bijection with the orbits of $\check
Q/(h+1)\check Q$ under $W$. 
\endproclaim
 
\demo
{Proof} Let $X= \{x\in V\mid (x,\alpha_i)\leq 1 \text{ for each } i\in
\{1,\dots ,n\} \text{ and }(x,\theta)\geq -2\}$ and
$\check\rho=\check\omega_1+\cdots+\check\omega _n$ be the half sum of positive
coroots. We have that $(\check\rho,\theta)=h-1$, thus $X$ is the simplex whose
vertices are $\check\rho$ and $\check\rho-(h+1)o_i$, for $i=1,\dots, n$. Hence
$X=t_{\check \rho} w_0(\overline C_{h+1})$. By Lemma 1 and the above Remark
there exists $w\in W_{af}$ such that $X=w(\overline C_{h+1})$. Such a $w$ gives
a bijection from $\overline C_{h+1}\cap \check Q$ to 
$D=X\cap \check Q$. If $\i\in \I$ and 
$w_\i= t_{\tau_\i} v_\i$, with $\tau_\i\in \check Q$ and $v_\i\in W$, then, using Proposition 3,
we obtain   that $w^{-1} v_\i^{-1}(\tau_\i)$ belongs to $\overline C_{h+1}\cap \check Q$ and 
$\i\mapsto w^{-1} v_\i^{-1}(\tau_\i)$ is a bijection between $\I$ and $\overline C_{h+1}\cap
\check Q$. This concludes the proof, since, as we observed in Section 2, the cosets of 
the elements in $\overline C_{h+1}\cap \check Q$ are a natural set of representatives 
of the orbits of $\check Q/(h+1)\check Q$ under the action of $W$.
\qed
\enddemo
\smallskip
We can explicitly determine the element $w$ 
which appears in the above proof.  
Indeed we shall compute $w^{-1}$.
If $\check \rho\in \check Q$, then trivially $w^{-1}=w_0 t_{-\check \rho}$. 
Otherwise, according to the proofs of Lemma 1  and of Theorem 1, there exists 
a unique $j\in J$ such that the vertex $\check \rho-(h+1)o_j=\check \rho-(h+1)\check\omega_j$
of $X$ belongs to $\check Q$. Now observe that $w_0^j$ maps $\{\alpha_i\mid i\not=j\}$ 
to $\{-\alpha_i\mid i\not=j\}$ and maps $\alpha_j$ and $\theta$ into positive roots. 
For any root $\alpha$ let $ht(\alpha)=(\alpha,\check\rho)$ be the height of $\alpha$. Then 
$ht(\theta)=h-1$, and, since $j\in J$, $ht(w_0^j(\theta))=ht(w^j_0(\alpha_j))-(h-2)$.
Since  
$w_0^j(\theta)$ is positive this implies that $ht(w^j_0(\alpha_j))=h-1$, hence 
$w^j_0(\alpha_j)=\theta$ and $w_0^j(\theta)=\alpha_j$. It is easily seen 
that this implies that $w_0^j t_{-\check\rho+(h+1)\check\omega_j}(X)=\overline C_{h+1}$. Hence 
in order to determine $w$ it suffices to determine the above $j$. 
\par 
Numbering the fundamental weights as in \cite{B}, by a direct computation we obtain: 
\medskip 
$A_n$: $\check \rho \in \check Q$ for $n$ even; $j={{n+1}\over 2}$ for $n$ odd.  
\par
$B_n$: $\check \rho \in \check Q$ for $n\equiv 0,3\ mod\  4$; $j=1$ for $n\equiv 1,2\ mod\  4$. 
\par
$C_n$: $j=n$.
\par
$D_n$: $\check \rho \in \check Q$  for $n\equiv 0,1\ mod\  4$; $j=1$ for $n\equiv 2,3\ mod\  4$.
\par
$E_7$: $j=7$.
\par
$E_6,\ E_8,\ F_4,\ G_2$:  $\check \rho \in \check Q$.

\medskip 
\heading 
\S4 The other bijections 
\endheading
\medskip

\noindent
{\bf (1)} {\sl A bijection between $ad$-nilpotent ideals  of $\bb$ and
antichains of the root poset}. 
In Section 3 we observed that any $ad$-nilpotent ideal of $\bb$ is a sum of
(positive) root spaces. 
For $\p\subseteq \Dp$,  set $\i_\p=\bigoplus\limits_{\a\in\p}\g_\a$. 
If $\i_\p$ is an ideal of $\bb$, then $\a\in\p,$ $\b\in\Dp,\,\a+\b\in\D$ implies $\a+\b\in\p$. 
If we endow $\Dp$ with the usual partial order  (i.e. $\a\leq \beta$ if
$\beta-\a=\sum_{\gamma\in\Dp}c_\gamma\gamma,\,c_\gamma\in\Bbb N$), then, by definition,  $\p$ is a
dual order ideal of $(\Dp,\leq)$.\par It is a general fact that, in a finite poset $P$,
dual order ideals and antichains (i.e. sets consisting of pairwise
non-comparable elements) are in canonical bijection: map the antichain
$\{a_1,\ldots,a_k\}$ to the dual  order ideal which is the union of the
principal dual order ideals $V_{a_1},\ldots,V_{a_k}$, where $V_a=\{b\in P\mid
b\geq a\}$; the inverse map sends a dual order ideal into the set of its minimal elements. 
It is clear that $\i_\p\mapsto \p$ is the required bijection.
\par\noindent{\bf Remark.}  In combinatorial literature the antichains of the
root poset $(\Dp,\leq)$ are called {\it nonnesting partitions} \cite{R, Remark
2}. This name derives from the analysis of the definition in type $A_n$. In that
case, write the positive roots with respect to the standard basis $\{\varepsilon_i\}_{i=1}^{n+1}$ of $\Bbb
R^{n+1}$, so that $\Dp=\{\varepsilon_i-\varepsilon_j\mid 1\leq i< j\leq n+1\}$.
Then to an antichain $A$ we can associate a partition of $\{1,\ldots,n+1\}$ by
putting in the same block $i,j$ whenever $\varepsilon_i-\varepsilon_j\in A$. It
turns out that partitions arising in this way  are the ones characterized by the following
property: if $a,e$ appear in a block $B$ and $b,d$ appear in a different block $B'$ where $a<b<d<e$, then there exists $c\in
B$ satisfying $b<c<d$.
\smallskip \noindent{\bf (2)} {\sl  A bijection between  antichains of the root
poset and $\oplus$-sign types of $\check\Delta$ or regions of the Catalan
arrangement  which are contained in the fundamental chamber}.
\par 
First we
recall the definition of $\oplus$-sign type for the  root system $\Delta$.
For $\a\in \Delta^+$ set $H_{\a,+}=\{v\in
V\mid (v,\check\a)>1\},\ H_{\a,0}=\{v\in V\mid 0<(v,\check\a)<1\},$ $
H_{\a,-}=\{v\in V\mid (v,\check\a)<0\}.$ Then a subset  $S\subset V$ is a
{\it sign type} (resp. $\oplus$-{\it sign type}) if it is of the form 
$S=\bigcap\limits_{\a\in\Dp}H_{\a,X_\a}$
for some collection $(X_\a)_{\a\in \Dp}$ with $X_\a\in\{+,0,-\}$ (resp. $X_\a\in\{+,0\}$).
\par  
We describe a bijection between dual order ideals and $\oplus$-sign types, 
according to Shi \cite{S, Theorem 1.4}. Given a dual
order ideal $\p\subseteq \Dp$, map it to   the $\oplus$-sign type $(X_{\check
\a})_{\check \a\in\check\Dp}$ defined by $$X_{\check \a}=\cases 0,\quad&\text{if
$\a\notin \p$}\\+,\quad&\text{if $\a\in \p$}.\endcases$$ 
This bijection appears also in a different context.  
Recall the two following remarkable arrangements
of real hyperplanes (cf. \cite{At2, \S3}). The  {\sl Shi arrangement} $\Cal S$, relative to $\D$,
is the set of hyperplanes of $V$ having equations $$(x,\a)=0,\qquad (x,\a)=1,\qquad \a\in\Dp;$$
the {\sl Catalan arrangement} $\Cal C$ is the set of hyperplanes of $V$ having equations
$$(x,\a)=0,\qquad (x,\a)=1,\qquad (x,\a)=-1,\qquad \a\in\Dp.$$ 
We call {\it regions} of the hyperplane arrangement the connected components 
of the complement in $V$ of the union of all hyperplanes in the arrangement.  
By the definition
of $\Cal S$ and $\Cal C$ it is clear that both arrangements have the same number
of regions inside the fundamental chamber of $W$.
A bijection between antichains in $\Dp$ and regions of 
$\Cal S$ or $\Cal C$ lying in the fundamental chamber 
(which in \cite{At1, 6.1} is attributed  to Postnikov)
can be made explicit mapping   an  antichain $A$ to the region 
$$X_A=\left\{x\in C_\infty\mid \cases (\b, x)>1\quad &\text{if $\b\geq\a$ 
for some $\a\in A$}\\ (\b,x)<1\quad &\text{otherwise}\endcases\right\}.$$

\smallskip 
We illustrate the above  bijections  in the easy case of a
root system of type $A_2$. $$ \alignat4 &\text{ideals in
$\I$}&&\text{antichains}\quad &&\text{regions of $\Cal C$ within $C_\infty$}\\ \\
&\i_1=0 &&\emptyset  && X_1=\{x\mid (x,\a_1)<1,(x,\a_2)<1,(x,\theta)<1\}\\
&\i_2=\g_\th &&\{\th\} && X_2=\{x\mid (x,\a_1)<1,(x,\a_2)<1,(x,\theta)>1\}\\
&\i_3=\g_{\a_1}\oplus\g_\th &&\{\a_1\}   &&X_3=\{x\mid
(x,\a_1)>1,(x,\a_2)<1,(x,\theta)>1\}\\ &\i_4=\g_{\a_2}\oplus\g_\th &&\{\a_2\} &&
X_4=\{x\mid (x,\a_1)<1,(x,\a_2)>1,(x,\theta)>1\}\\
&\i_5=\g_{\a_1}\oplus\g_{\a_2}\oplus\g_\th\quad &&\{\a_1,\a_2\}\quad\quad &&
X_5=\{x\mid (x,\a_1)>1,(x,\a_2)>1,(x,\theta)>1\}\endalignat$$
\smallskip
\noindent
{\bf Remark.} It is worthwhile to recall that formula (1) 
also counts the number of conjugacy classes of elements of order
dividing $h+1$ in a maximal torus $T$ of the connected simply connected simple algebraic group $G$ corresponding to $\g$.
Indeed these classes are in  bijection with  $W$-orbits on  $\check Q/r\check Q$. In fact, 
regard coroots as cocharacters of $T$, i.e. as morphism
of algebraic groups  $\Bbb C^*\to T$. Fix a primitive $r$-th root of unity 
$z$; then, given $\tau\in \check Q$, the map $\tau\mapsto\tau(z)$ is
bijection from $\check Q/r\check Q$ to $T_r=\{t\in T\mid t^r=1\}$ and induces a 
bijection between the $W$-orbits in $\check Q/r\check Q$ and the conjugacy classes of 
elements in $T_r$.
\medskip 

\medskip 
\heading 
\S5 Examples
\endheading

\medskip 
We illustrate the bijection of Theorem 1 when $\D$ is of type  $A_2$ or $B_2$.
For this purpose we first need to give explicitly the elements $w_\i\in \widehat W$ 
corresponding to the ideals $\i\in \I$.  
\par 
In the case of $A_2$ the map from $\I$ into $\widehat W$ is given as follows. 
$$
\alignat3 &\text{ideals in $\I$}&&N(w_\i)\quad &&w_\i\\ \\ &\i_1=0 &&\emptyset 
&& 1\\ &\i_2=\g_\th &&\{-\th+\d\} && s_0\\ &\i_3=\g_{\a_1}\oplus\g_\th
&&\{-\th+\d,-\a_1+\d\}   &&s_0s_2\\ &\i_4=\g_{\a_2}\oplus\g_\th
&&\{-\th+\d,-\a_2+\d\} && s_0s_1\\
&\i_5=\g_{\a_1}\oplus\g_{\a_2}\oplus\g_\th\quad
&&\{-\th+\d,-\a_2+\d,-\th+2\d,-\a_1+\d\}\quad\quad &&s_0s_1s_2s_1
\endalignat
$$
We have $h=3$ and $\check Q=Q=\Bbb Z\a_1\oplus\Bbb Z \a_2$; we
have also $\check \rho=\rho =\th=\a_1+\a_2$. The bijections of Proposition 3 and 
Theorem 1 look as follows (regarding the bijection of Theorem 1, we write down 
the element of $\overline C_{h+1}\cap \check Q$ corresponding to each ideal). 

$$
\alignat3 &w_\i=t_{\tau_\i}v_\i\ &&v_\i^{-1}(\tau_\i) &&w_0t_{-\check\rho}(v_\i^{-1}(\tau_\i)) \\ \\ 
&1 &&0\qquad\qquad &&\th \\ &s_0=t_\th s_1s_2s_1\qquad\qquad 
&&-\th\qquad\qquad &&
2\th\\ &s_0s_2=t_\th s_2s_1 &&-\a_1\qquad\qquad && 
\a_1+2\a_2\\ &s_0s_1=t_\th s_1s_2 &&-\a_2\qquad\qquad  &&
2\a_1+\a_2\\ &s_0s_1s_2s_1=t_\th &&\th \qquad\qquad &&0 
\endalignat
$$

Now  we consider the root type $B_2$. Here $h=4$, and  $\check
Q=\Bbb Z\a_1+2\Bbb Z\a_2$. Then  $\check \rho=2\a_1+3\a_2\notin\check Q$, but  
$\check\omega_1=\a_1+\a_2$ so that
$\check\rho-5 \check\omega_1=-3 \a_1- 2\a_2 \in \check Q$. Moreover, in the notation of 
the proof of Theorem 1,  $w_0^{j}=w_0^{1}=s_2$. 
The injection of $\I$ in $\widehat W$ is given as follows (for
shortness we don't write $N(w_\i)$).

$$
\alignat2 &\i_1=0\qquad
\qquad\qquad\qquad\qquad\qquad&&w_{\i_1}=1\\ &\i_2=\g_{\th}\quad
&&w_{\i_2}=s_0\\ &\i_3=\g_{\a_1+\a_2}\oplus\g_{\th}\quad &&w_{\i_3}=s_0s_2\\
&\i_4=\g_{\a_1}\oplus\g_{\a_1+\a_2}\oplus\g_{\th}\quad &&w_{\i_4}=s_0s_2s_0\\
&\i_5=\g_{\a_2}\oplus\g_{\a_1+\a_2}\oplus\g_{\th}\quad &&w_{\i_5}=s_0s_2s_1s_2\\
&\i_6=\frak n\quad &&w_{\i_6}=s_0s_2s_1s_2s_0s_2s_0
\endalignat
$$ 

The bijection with $\overline C_{h+1}\cap \check Q$ is obtained in the following table.
$$\alignat3
&w_\i=t_{\tau_\i}v_\i &&v_\i^{-1}(\tau_\i) \quad
&&w_0^{j}t_{-\check \rho+(h+1)\check \omega_{j}}(v_\i^{-1}(\tau_\i))\\ \\
&w_{\i_1}=1\qquad\qquad &&0\qquad\qquad &&3\check\a_1+2\check \a_2 \\ &w_{\i_2}=t_{\check
\th}s_2s_1s_2\qquad\qquad &&-\check\a_1-\check\a_2\qquad\qquad && 2\check\a_1+2\check\a_2\\
&w_{\i_3}=t_{\check \th}s_2s_1\qquad\qquad &&-\check\a_1\qquad\qquad && 
2\check\a_1+\check\a_2\\ &w_{\i_4}=t_{2\check\a_1+\check\a_2}s_1s_2s_1 \qquad\qquad&&
-2\check\a_1-\check\a_2\qquad\qquad && \check\a_1+\check\a_2\\ &w_{\i_5}=t_{\check
\th} \qquad\qquad&&\check\a_1+\check\a_2\qquad\qquad  && 4\check\a_1+2\check\a_2\\
&w_{\i_6}=t_{3\check\a_1+2\check\a_2}s_1s_2s_1\qquad\qquad &&-3\check\a_1-\check\a_2
\qquad\qquad &&0 
\endalignat
$$

\medskip 
\par 
\newpage
\Refs \widestnumber\key {PPPPP} 
\bigskip \bigskip 


\ref\key{\bf At1}\by C. A. Athanasiadis \yr 1998 \paper
On Noncrossing and nonnesting partitions for classical reflection groups\jour 
Electronic J. of Comb\. \vol 5 \pages \#R42\endref 

\ref\key{\bf At2}\by C. A.
Athanasiadis \paper Deformations of Coxeter hyperplane arrangements and their
characteristic polynomials\inbook Arrangements --- Tokyo 1998\eds M. Falk and H.
Terao\publ Advanced Studies in Pure Mathematics, in press\endref 

\ref\key{\bf
B}\by N. Bourbaki \book Groupes et algebres de Lie \publ Hermann\publaddr
Paris\yr 1968 \endref 

\ref\key{\bf CP1}\by P. Cellini and P. Papi  \paper
$ad$-nilpotent ideals of a Borel subalgebra\jour J. Algebra \yr 2000\vol
225\pages130--141 \endref 

\ref\key{\bf CP2}\by P. Cellini and P. Papi  \paper
$ad$-nilpotent ideals of a Borel subalgebra III\jour in preparation  \endref 

\ref\key{\bf Ha}\by M. Haiman\paper Conjectures on the
quotient ring by diagonal invariants\jour J. Algebraic Combin. \vol 3 \yr
1994\pages 17--76  \endref 

\ref
\key{\bf Hu}\by J. E. Humphreys \yr 1990 \book
Reflection groups and Coxeter groups\publ Cambridge University Press\publaddr
Cambridge
\endref 

\ref
\key{\bf Kac}\by V.G. Kac \book Infinite Dimensional Lie
Algebras \publ Cambridge University Press\yr 1990 
\endref 

\ref\key{\bf IM}\by N. Iwahori, H. Matsumoto\pages5--48
\paper On some Bruhat decomposition and the structure of the Hecke rings of p-adic
Chevalley groups \yr1965\vol 25
\jour Inst. Hautes \'Etudes Sci. Publ. Math.
\endref

\ref
\key{\bf Ko1}\by
B. Kostant \paper Eigenvalues of a Laplacian and commutative Lie subalgebras
\yr1965\vol 3, suppl. 2 \jour Topology\pages 147--159
\endref 

\ref
\key{\bf Ko2}\by B. Kostant \paper The Set of
Abelian ideals of a Borel Subalgebra, Cartan Decompositions, and Discrete Series
Representations
\yr1998\vol5 
\jour Internat. Math. Res. Notices\pages 225--252\endref

\ref
\key{\bf KOP}\by  C. Krattenthaler, L. Orsina and P. Papi  \paper enumeration of
$ad$-nilpotent $\frak b$-ideals for simple Lie algebras \jour Adv\. Appl\.
Math\.  (to appear)\finalinfo {\tt math.RA/0011023} 
\endref 
\ref
\key{\bf PR}\by D. Panyushev, G. R\"ohrle  
\paper Spherical orbits and abelian ideals \jour 
Adv\.  Math\. 
\vol 159\pages 229--246\yr 2001 
\endref 


\ref
\key{\bf R}\by V.    Reiner \yr 1997 \paper
Non-crossing partitions for classical reflection groups\jour Discrete Math\.\vol
177\pages 195--222
\endref 

\ref
\key{\bf S}\by J. Shi \paper The number of
$\oplus$-sign types \jour Quart. J. Math. Oxford\vol 48\yr
1997\pages 93--105
\endref 

\ref
\key{\bf So}\by E. Sommers \yr 1997\jour Transform.
Groups \paper A family of affine Weyl group representations\vol 2\pages
375--390
\endref


\endRefs
\enddocument